\documentclass[12pt]{article}
\usepackage{amsmath,amssymb,amsfonts,amsthm,graphicx,float,fancyhdr}
\usepackage[a4paper,left=3cm,right=3cm,top=3cm,bottom=3cm]{geometry}
\usepackage{hyperref}
\usepackage[pagewise]{lineno}
\newtheorem{theorem}{Theorem}[section]

\newtheorem{conjecture}[theorem]{Conjecture}
\numberwithin{equation}{section}
\newenvironment{dedication}
        {\vspace{3ex}\begin{quotation}\begin{center}\begin{em}}
        {\par\end{em}\end{center}\end{quotation}}

\title{ A Lochs-Type Approach via Entropy in Comparing the Efficiency of Different Continued Fraction Algorithms}
\author{
    Dan Lascu\footnote{e-mail: lascudan@gmail.com, corresponding author}\\
    \emph{\small Mircea cel Batran Naval Academy, 1 Fulgerului, 900218 Constanta,
    Romania}
    and\\
    Gabriela Ileana Sebe\footnote{e-mail: igsebe@yahoo.com.} \\
    \emph{\small Politehnica University of Bucharest, Faculty of Applied Sciences},\\
    \emph{\small Splaiul Independentei 313, 060042, Bucharest, Romania} and \\
    \emph{\small Gheorghe Mihoc-Caius Iacob Institute of Mathematical Statistics } \\
    \emph{\small and Applied Mathematics}, \\
    \emph{\small Calea 13 Sept. 13, 050711 Bucharest, Romania}
    }
\begin{document}
\maketitle

\begin{dedication}
{Dedicated to the memory of our mentor, Professor Marius Iosifescu.}
\end{dedication}

\abstract{We investigate the efficiency of several types of continued fraction expansions of a number in the unit interval using a generalization of Lochs theorem from 1964. Thus, we aimed to compare the efficiency by describing the rate at which the digits of one number-theoretic expansion determine those of another. We study Chan's continued fractions, $\theta$-expansions, $N$-continued fractions and R\'enyi-type continued fractions. A central role in fulfilling our goal is played by the entropy of the absolutely continuous invariant probability measures of the associated dynamical systems.}
\section{Introduction}\label{sec1}
The purpose of this paper is to compare the efficiency of some continued fractions in approximating a real number in the unit interval.
For any irrational $x \in (0,1)$ suppose we have two known expansions $x=[a_1, a_2, \ldots]_1$ and $x=[b_1, b_2, \ldots]_2$.
A natural question is: which of these continued fraction expansions is more ``efficient''?
In mathematical terms ``efficiency'' means which of the two sequences $[a_1, a_2, \ldots, a_n]_1$, respectively $[b_1, b_2, \ldots, b_n]_2$ converges faster to $x$ as $n \rightarrow \infty$?
It is therefore relevant to ask how much information (e.g., in terms of the digits in the second expansion) can be determined once we know $n$ digits of the first expansion.
Suppose we approximate $x$ by keeping the first $n$ digits of its first expansion.
In order to attain the same degree of accuracy, we need  to keep the first $m=m(n,x)$ digits of the second expansion.
What can we say about the ratio $m(n,x)/n$ in general?
We shall see that the relative speed of approximation of two different expansions (almost everywhere) is related to the quotient of the entropies of the transformations that generate these expansions.
The strategy used in this paper fits for all pairs of number-theoretic fibred maps (NTFMs) for which the generating partition has finite entropy.

\subsection{Motivation} \label{sec1.1}
The ancient charms of number theory have been replaced by the modern fascination of an algorithmic thinking.
When analytic number theory was being formulated in the nineteenth century, probability theory was not yet a reputable branch of mathematics.
Nowadays, probabilistic techniques are routinely used in number theory.
The field of probabilistic number theory is currently evolving very rapidly and uses more and more refined probabilistic tools and results.
It is well-known that continued fractions lie at the heart of a number of classical algorithms like Euclid's greatest common divisor or the lattice reduction algorithm of Gauss \cite{IK-2002}. So, continued fractions arise naturally in the theory of approximation to real numbers by rationals.
To this day the Gauss map, on which metrical theory of regular continued fraction is based, has fascinated researchers from various branches of mathematics and science, with many applications in computer science, cosmology and chaos theory. Apart from the regular continued fraction expansion, very many other continued fraction expansions were studied. In the last century, mathematicians broke new ground in this area.
Since there are several continued fraction algorithms, we ask ourselves which of them provides the best approximation of a real number.
The representation of a real number by a continued fraction can be viewed as a source of information about the number.
For this purpose we need the notion of entropy, a rigorous tool at the crossroads between probability, information theory and dynamical systems \cite{Ebrahimzadeh,Ryabko-2019}.

\subsection{Entropy} \label{sec1.2}
As is well known, entropy is an important concept of information in physics, chemistry, and information theory \cite{PY-1998}.
The connection between entropy and the transmission of information was first studied by Shannon in \cite{Shannon-1948}.
Thus, the entropy can be seen as a measure of randomness of the system, or the average information acquired under a single application of the underlying map.
Entropy also plays an important role in ergodic theory.
Thus Shannon's probabilistic notion of entropy was first introduced into the ergodic theory by Kolmogorov \cite{Kolmogorov-1958} via a measure theoretic approach. The contribution of Kolmogorov to modern dynamics was the discovery of the concept of entropy, which was made rigorous by Sinai \cite{Sinai-1959}.
This concept provides an important generalization of Shannon entropy. Kolmogorov-Sinai (K-S) entropy measures the maximal loss of information for the iteration of finite partitions in a measure preserving transformation.
The concept has shown its strength through the adequate answers to problems in the classification of dynamical systems.
Two metrically isomorphic dynamical systems have the same K-S entropy, so this concept is a tool for distinguishing non-isomorphic(non-conjugate) dynamical systems.

We briefly review this very important concept of K-S entropy in Ergodic Theory.
Given a measure preserving system $\left( X, \mathcal{X}, \mu, T \right)$,
we say $\alpha = \{ A_i : i \in I \}$ is a \textit{partition} of $X$ if
$X=\bigcup_{i \in I} A_i$, where $A_i \in \mathcal{X}$ for each $i \in I$ and $A_i \cap A_j = \emptyset$ for $i \neq j$, $i, j \in I$.
Here $I$ is a finite or countable index set.
For a partition $\alpha$ of $X$, we define \textit{the entropy of the partition} $\alpha$ as
\begin{equation}\label{1.01}
H (\alpha) := - \sum_{A \in \alpha} \mu (A) \log \mu (A).
\end{equation}
In this definition $T$ does not appear.
But, the entropy of the dynamical system is defined by the entropy of the transformation $T$ as follows.

Given a partition $\alpha$, we consider the partition
\begin{equation}\label{1.02}
\alpha_n := \bigvee^{n-1}_{i=0} T^{-i}\alpha = \left\{ \bigcap^{n-1}_{i=0} T^{-i}A_i : A_i \in \alpha, i=0, 1, \ldots, n-1 \right\}.
\end{equation}
Then \textit{the entropy of transformation} $T$ w.r.t. $\alpha$ is given by
\begin{equation}\label{1.03}
  h(\alpha, T) := \lim_{n \rightarrow \infty} \frac{1}{n} H(\alpha_n).
\end{equation}
The entropy of $T$ is defined as
\begin{equation}\label{1.003}
  h(T) := \sup\left\{h(\alpha, T): \alpha \mbox{ partition of } X \mbox{ such that } H(\alpha)< \infty\right\}.
\end{equation}
In general, it does not seem possible to calculate the entropy straight from its definition.
First, let us define that a partition $\alpha$ of $X$ is a generator w.r.t. a non-invertible transformation $T$ if
\begin{equation}\label{1.003}
\sigma \left( \bigvee^{\infty}_{i=0} T^{-i}\alpha \right) = \mathcal{X} \mbox{ up to sets of } \mu\mbox{-measure zero}.
\end{equation}
In computation of $h(T)$ the K-S Theorem is very useful.
For completeness, we recall this theorem.
Let $\alpha$ be a partition of $X$ such that $H(\alpha)< \infty$.
If $\alpha$ is a generator w.r.t. $T$, then $h(T)=h(\alpha,T)$.

We also have the following classical Shannon-McMillan-Breiman Theorem \cite{DK-2001}.
Let $\left( X, \mathcal{X}, \mu, T \right)$ be an ergodic measure preserving system and let $\alpha$ be a finite or countable partition of $X$ satisfying $H(\alpha)< \infty$.
The Shannon-McMillan-Breiman theorem says if $A_n(x)$ denote the unique element $A_n \in \alpha_n$ such that $x \in A_n$,
then for almost every $x \in X$ we have:
\begin{equation}\label{1.04}
  \lim_{n \rightarrow \infty} -\frac{1}{n} \log \mu (A_n(x)) = h(\alpha, T).
\end{equation}
In 1964 Rohlin \cite{Rohlin-1964} showed that, when R\'enyi's condition is satisfied the entropy of a $\mu$-measure preserving operator $T : X \rightarrow X$ is given by the formula:
\begin{equation} \label{1.1}
h(T) := \int_{X} \log\left|T'(x)\right| \mathrm{d}\mu(x).
\end{equation}
The \textit{R\'eny's condition} means that there is a constant $C \geq 1$ such that
\begin{equation}\label{1.05}
  \sup_{x,y \in T^n(\alpha_n)} \frac{\left| u'_n(x) \right|}{\left| u'_n(y) \right|} \leq C,
\end{equation}
where $ u_n := \left( \left.T^n\right|_{\alpha_n} \right)^{-1}$.
%

\section{Lochs' Theorem} \label{sec2}
Expansions that furnish increasingly good approximations to real numbers are usually related to dynamical systems.
In 1964, G. Lochs \cite{Lochs-64} compared the decimal expansion and the regular continued fraction (RCF) expansion.
Although comparing dynamical systems seems difficult, using detailed knowledge of the continued fraction operator, Lochs was able to relate the relative speed of approximation of decimal and regular continued fraction expansions (almost everywhere) to the quotient of the entropies of their dynamical systems.
Thus, roughly $97$ RCF digits are determined by $100$ decimal digits which indicates
that the RCF expansion is slightly more efficient compared to the decimal expansion at representing irrational numbers.
\subsection{Decimal expansions} \label{sec2.1}
It is a common known fact that, any real number $x \in [0, 1)$ can be written as
\begin{equation} \label{2.1}
x = \sum_{i=1}^{\infty} \frac{d_i(x)}{10^i},
\end{equation}
where $d_i=d_i(x) \in \{0, 1, 2, \ldots, 9\}$ for $i \in \mathbb{N}_+:=\{1, 2, 3, \ldots \}$.
The representation of $x$ in (\ref{2.1}), denoted by $x=0.d_1d_2\ldots$ is called the \textit{decimal expansion} of $x$.
We can generate the decimal expansions by iterating the \textit{decimal map}
\begin{equation} \label{2.2}
T_d: [0, 1) \rightarrow [0, 1); \quad T_d(x) = 10x - \lfloor 10x \rfloor,
\end{equation}
where $\lfloor \cdot \rfloor$ denotes the floor (or entire) function.
In other words, $T_d$ is given by
\begin{equation} \label{2.3}
T_d(x) = 10x - i \, \, \mbox{ if } \,  \frac{1}{10} \leq x < \frac{i+1}{10}, \, i=0, 1, 2, \ldots, 9.
\end{equation}
Thus, we obtain:
\begin{equation} \label{2.4}
x = \frac{d_1}{10} + \frac{d_2}{10^2} + \ldots + \frac{d_n}{10^n} + \frac{T_d^n(x)}{10^n},
\end{equation}
where
\begin{equation} \label{2.5}
d_1=d_1(x)= \lfloor 10x \rfloor
\end{equation}
and
\begin{equation} \label{2.6}
d_n=d_n(x)=d_1\left(T^{n-1}_d(x)\right), \quad n \geq 2.
\end{equation}
Since $ 0 \leq T_d^n(x) <1$, we obtain
\begin{equation} \label{2.7}
\sum_{i=1}^{n} \frac{d_i}{10^i} \rightarrow x \quad \mbox{as } n \rightarrow \infty.
\end{equation}
\subsection{Regular continued fraction expansions} \label{sec2.2}
Beside decimal expansions there are many more possible representations of real numbers in terms of a sequence of integers. We refer to the regular continued fractions as another famous example. Any irrational number $x \in [0, 1)$ has a unique \textit{regular continued fraction expansion}
\begin{equation} \label{2.8}
x=\displaystyle \frac{1}{a_1+ \displaystyle \frac{1}{a_2+ \ddots}},
\end{equation}
where $a_n \in \mathbb{N}_+$ for any $n \geq 1$.
This expansion is obtained by applying repeatedly the \textit{Gauss map} or the \textit{regular continued fraction transformation}
\begin{equation} \label{2.9}
T_G: [0, 1) \rightarrow [0, 1); \quad T_G(0)=0, \, T_G(x) = \frac{1}{x} - \left\lfloor \frac{1}{x} \right\rfloor.
\end{equation}
Therefore, it follows that the digits $a_1, a_2, \ldots$ are related by
\begin{equation} \label{2.10}
a_n=a_n(x)=a_1\left(T^{n-1}_G(x)\right), \quad n \geq 2,
\end{equation}
where
\begin{equation} \label{2.11}
a_1=a_1(x)= \left\lfloor \frac{1}{x} \right\rfloor.
\end{equation}
If we denote by $[a_1, a_2, \ldots]_G$ the expansion in (\ref{2.8}), we have
\begin{equation} \label{2.12}
[a_1, a_2, \ldots, a_n]_G \rightarrow x \quad \mbox{as } n \rightarrow \infty.
\end{equation}
\subsection{Comparing the efficiency of decimal expansion and RCF expansion} \label{sec2.3}
In \cite{Lochs-64}, Lochs answered the question which of these developments is more efficient, namely,
which of the two sequences in (\ref{2.7}) and (\ref{2.12}) converges faster to $x$ as $n \rightarrow \infty$.

Suppose that the irrational number $x \in (0, 1)$ has the decimal expansion $x = 0.d_1d_2\ldots$
and the RCF expansion $x = [a_1, a_2, \ldots]_G$.
Let $y=0.d_1d_2 \ldots d_n$ be the rational number determined by the first $n$ decimal digits of $x$, and let $z = y+10^{-n}$.
Then, $[y, z]$ is the $n$-th order decimal cylinder containing $x$, which we also denote by $C_n(x)$.
Now let $y=[b_1, \ldots, b_l]_G$ and $z=[c_1, c_2, \ldots, c_k]_G$ be the RCF expansions of $y$ and $z$.
Let
\begin{equation} \label{2.13}
m(n, x) := \max\{ i \leq \max(l, k): \mbox{ for all } j\leq i, \, b_j=c_j  \}.
\end{equation}
In other words, $m(n,x)$ is the largest integer such that $C_n(x) \subset D_{m(n,x)}(x)$, where $D_j(x)$ denotes the $j$-th order RCF cylinder containing $x$.
Lochs \cite{Lochs-64} proved that, for almost every irrational $x \in (0, 1)$, we have
\begin{equation} \label{2.14}
\lim_{n \rightarrow \infty} \frac{m(n, x)}{n} = \frac{6 \log 2 \log 10}{\pi^2} = 0.97027014\ldots
\end{equation}
For example, Lochs \cite{Lochs-63} has computed the first $968$ RCF digits of $\pi$ from its first $1000$ decimals,
and Brent and McMillan \cite{Brent} have computed the first $29200$ RCF digits of Euler–Mascheroni constant from its first $30100$ decimals.
\subsection{Extended Lochs' theorem}\label{sec2.4}
Dajani and Fieldsteel \cite{DF-2001} proved a generalization of Lochs' theorem showing that
we can compare any two expansions of numbers which are generated by number-theoretic fibred maps, i.e., surjective interval maps $S:[0, 1) \rightarrow [0, 1)$ that satisfy the following conditions:
\begin{itemize}
  \item[(1)] there exists a finite or countable partition of $[0, 1)$ into intervals such that $S$ restricted to each interval is strictly monotonic and continuous,

    \noindent and

  \item[(2)] $S$ is ergodic with respect to Lebesgue measure $\lambda$, and there exists an $S$ invariant probability measure $\mu$ equivalent to $\lambda$ (i.e., $\mu(A) = 0$ if and only if $\lambda(A)=0$ for all Lebesgue sets $A$) with bounded density (both $\displaystyle\frac{\mathrm{d}\mu}{\mathrm{d}\lambda}$ and $\displaystyle\frac{\mathrm{d}\lambda}{\mathrm{d}\mu}$ are bounded).
\end{itemize}
Consider NTFMs $S_1$ and $S_2$ on $[0,1)$ with invariant measures $\mu_1$ and $\mu_2$ (equivalent to Lebesgue measure) and with partitions $P$ and $Q$, respectively. Denote by $P_n(x)$ the $n$-th order cylinder that contains $x \in [0,1)$ (a similar definition for $Q_n(x)$).
Let
\begin{equation}
m_{S_1,S_2} (n,x)=\sup \{ m: P_n(x) \subset Q_m(x) \}.
\end{equation}
Under the conditions just stated and with the understood that $m_{S_1,S_2} (n,x)$ is exactly the number of digits in the $S_2$-expansion of $x$ that can be determined from knowing the first $n$ digits in the $S_1$-expansion, we have
\begin{equation} \label{2.15}
\lim_{n \rightarrow \infty} \frac{m_{S_1,S_2}(n, x)}{n} = \frac{h(S_1)}{h(S_2)} \quad \lambda - \mathrm{a.e.}
\end{equation}
where $h(S_1)$ and $h(S_2)$ denote the entropy of $S_1$ and $S_2$, respectively, with $h(S_1)>0$ and $h(S_2)>0$.

\section{Other continued fraction expansions} \label{sec3}
Apart from the RCF expansions there is a wide variety of continued fraction expansions.
Here we mention only a few of the expansions studied from the metrical point of view by the authors over time, namely,
Chan's continued fractions (in \cite{Lascu-2013,Sebe-2010}),
$\theta-$expansions (in \cite{SebeLascu-2014,Sebe-2017}),
$N-$continued fraction expansions (in \cite{Lascu-2016,SL-2020-4}),
R\'enyi-type continued fraction expansions (in \cite{LS-2020-1,LS-2020-2,SL-2020-3}).

In this paper, we ask the question which of these expansions is more effective.
In order to apply the result of Dajani and Fieldsteel \cite{DF-2001} from section \ref{sec2.4}, we briefly present these expansions, calculating the entropy of each map that generates these expansions.
\subsection{Chan's continued fractions} \label{sec3.1}
Fix an integer $\ell \geq 2$. Then any $x \in \left[0, 1\right)$ can be written in the form
\begin{equation} \label{3.1}
x = \frac{\ell^{-a_1(x)}}{\displaystyle 1+\frac{(\ell-1)\ell^{-a_2(x)}}{\displaystyle 1 + \frac{(\ell-1)\ell^{-a_3(x)}}{\displaystyle 1 + \ddots} }}
=:[a_1(x), a_2(x), a_3(x), \ldots]_{\ell},
\end{equation}
where $a_n(x)$'s are non-negative integers.
Define $ \frac{p_{\ell,n}(x)}{q_{\ell,n}(x)} := [a_1(x), a_2(x), \ldots, a_n(x)]_{\ell}$ the $n$-convergent of $x$ by truncating the expansion on the right-hand side of (\ref{3.1}), that is,
\begin{equation} \label{3.2}
[a_1(x), a_2(x), \ldots, a_n(x)]_{\ell} \to x \quad (n \rightarrow \infty).
\end{equation}
This continued fraction is associated with the following transformation $T_{\ell}$
on $[0,1]$:
\begin{equation}
T_{\ell} (x) := \left\{\begin{array}{lll}
\displaystyle \frac{\ell^{ \frac{\log x^{-1}}{\log \ell} - \left\lfloor\frac{\log x^{-1}}{\log \ell}\right\rfloor} - 1}{\ell-1}, & \hbox{if} & x \neq 0 \\
0, & \hbox{if} & x = 0.
\end{array} \right. \label{3.3}
\end{equation}
We notice that $T_{\ell}$ maps the set of irrationals in $[0, 1]$ into itself.
For any $x \in (0,1)$ put
\begin{equation} \label{3.4}
a_n=a_n(x) = a_1\left(T_{\ell}^{n-1}(x)\right), \quad n \in {\mathbb{N}}_+,
\end{equation}
with $T_{\ell}^0 (x) = x$ and
\begin{equation} \label{3.5}
a_1=a_1(x) = \left\{\begin{array}{lll}
\lfloor\log x^{-1} / \log \ell \rfloor, & \hbox{if} & x \neq 0 \\
\infty, & \hbox{if} & x = 0.
\end{array} \right.
\end{equation}

The transformation $T_{\ell}$ which generates the continued fraction expansion (\ref{3.1}) is ergodic with respect to an invariant probability measure, $G_{\ell}$, where
\begin{equation}\label{3.6}
 G_{\ell} (A) := k_{\ell} \int_{A} \frac{\mathrm{d}x}{((\ell-1)x+1)((\ell-1)x+\ell)}, \quad A \in {\mathcal{B}}_{[0, 1]},
\end{equation}
with
\begin{equation}\label{3.7}
k_{\ell} := \frac{(\ell-1)^2}{\log \left(\ell^2/(2\ell-1)\right)}
\end{equation}
and $\mathcal{B}_{[0, 1]}$ is the $\sigma$-algebra of Borel subsets of $[0, 1]$.

An $n$-block $(a_1, a_2, \ldots, a_n)$ is said to be \textit{admissible} for the expansion in (\ref{3.1}) if there exists
$x \in [0, 1)$ such that $a_i(x)=a_i$ for all $1 \leq i \leq n$.
If $(a_1, a_2, \ldots, a_n)$ is an admissible sequence, we call the set
\begin{equation} \label{3.06}
I_{\ell} (a_1, a_2, \ldots, a_n) = \{x \in [0, 1]:  a_1(x) = a_1, a_2(x) = a_2, \ldots, a_n(x) = a_n \},
\end{equation}
\textit{the $n$-th order cylinder}.

Define $(u_{\ell,i})_{i \in \mathbb{N}}$ by
\begin{equation}\label{3.007}
 u_{\ell,i}: [0, 1] \rightarrow [0, 1]; \quad   u_{\ell,i}(x) := \frac{\ell^{-i}}{1+(\ell-1)x}.
\end{equation}
For each $i \in \mathbb{N}$, $u_{\ell,i} = \left( \left.T_{\ell}\right|_{I_{\ell}(a_1)} \right)^{-1}$.
Let
\begin{equation}\label{3.008}
 u_{\ell,a_1a_2 \ldots a_n}(t) := \left(u_{\ell,a_1} \circ u_{\ell,a_2} \circ \ldots \circ u_{\ell,a_n}\right)(t) =
 \frac{\ell^{-a_1(x)}}{\displaystyle 1+\frac{(\ell-1)\ell^{-a_2(x)}}{\displaystyle 1 + \ddots +\frac{(\ell-1)\ell^{-a_n}}{1+(\ell-1)t}}}.
\end{equation}
We observe that $u_{\ell,a_1a_2 \ldots a_n}=\left( \left.T^n_{\ell}\right|_{I_{\ell}(a_1, a_2, \ldots, a_n)} \right)^{-1}$.
Therefore
\begin{equation} \label{3.009}
I_{\ell} (a_1, a_2, \ldots, a_n) = \{ u_{\ell,a_1 \ldots a_n}(t) : t \in [0, 1) \}
\end{equation}
which is an interval with the endpoints
$\frac{p_{\ell,n}}{q_{\ell,n}}$ and $\frac{p_{\ell,n}+(\ell-1)\ell^{a_n}p_{\ell,n-1}}{q_{\ell,n} + (\ell-1)\ell^{a_n}q_{\ell,n-1}}$.
Such intervals form a partition of $[0, 1]$.

Before applying Rohlin's formula, we must check R\'enyi condition (\ref{1.05}).
We will use directly, without mentioning them here, some properties proved in \cite{Lascu-2013}.
Thus, we have
\begin{equation}\label{3.0010}
\frac{\left| u'_{\ell,a_1 \ldots a_n}(t) \right|}{\left| u'_{\ell,a_1 \ldots a_n}(r) \right|} =
\left( \frac{q_{\ell,n}+r(\ell-1)\ell^{a_n}q_{\ell,n-1}}{q_{\ell,n}+t(\ell-1)\ell^{a_n}q_{\ell,n-1}} \right)^2 \leq
\left( \frac{q_{\ell,n}+(\ell-1)\ell^{a_n}q_{\ell,n-1}}{q_{\ell,n}} \right)^2 \leq \ell^2.
\end{equation}

Applying Rohlin's formula (\ref{1.1}) we obtain:
\begin{equation}\label{3.8}
\begin{split}
   h(T_{\ell})  =& \int_{0}^{1} \log \left|T'_{\ell}(x)\right| \mathrm{d}G_{\ell} = \int_{0}^{1} \log \left( \frac{\ell^{-a_1(x)}}{(\ell-1)x^2} \right)\mathrm{d}G_{\ell} \\
           =& \int_{0}^{1} \left( 2 \log(1/x) - a_1(x) \log \ell - \log(\ell-1)\right) \mathrm{d}G_{\ell} \\
           =& 2 k_{\ell} \int_{0}^{1} \frac{\log(1/x)}{((\ell-1)x+1)((\ell-1)x+\ell)} \mathrm{d}x \\
            & - k_{\ell} (\log \ell) \int_{0}^{1} \frac{a_1(x)}{((\ell-1)x+1)((\ell-1)x+\ell)}\mathrm{d}x - \log(\ell-1).
\end{split}
\end{equation}
As examples, we have the values of $h(T_{\ell})$ for different values of ${\ell}$ and the graph of $h(T_{\ell})$ in Table~\ref{tab1}.

\subsection{$\theta-$expansions} \label{sec3.2}
For a fixed irrational $\theta \in (0,1)$, we consider a generalization of the Gauss map,
$T_{\theta}: [0,\theta] \to [0,\theta]$ defined as
\begin{equation} \label{3.9}
T_{\theta}(x):=
\left\{
\begin{array}{ll}
{\displaystyle \frac{1}{x} - \theta \left \lfloor \frac{1}{x \theta} \right\rfloor}, &
{\displaystyle \mbox{if } x \in (0, \theta],}\\
0, & \mbox{if } x=0.
\end{array}
\right.
\end{equation}
The transformation $T_{\theta}$ is connected with the $\theta-$expansion for a number in $(0, \theta)$ as follows.
The numbers $\theta \left\lfloor \frac{1}{y \theta}  \right\rfloor$ obtained by taking $y$ successively equal to
$x$, $T_{\theta}(x)$, $T_{\theta}^2(x)$, \ldots, lead to the $\theta$-expansion of $x$ as
\begin{equation}
x = \frac{1}{\displaystyle \vartheta_1\theta
+\frac{1}{\displaystyle \vartheta_2\theta
+ \frac{1}{\displaystyle \vartheta_3\theta + \ddots} }} =: [\vartheta_1, \vartheta_2, \vartheta_3, \ldots]_{\theta}, \label{3.10}
\end{equation}
where $\vartheta_n \in \mathbb{N}_+$ for all $n \in \mathbb{N}_+$.
The positive integers $\vartheta_n=\vartheta_n(x)=\vartheta_1\left(T_{\theta}^{n-1}(x)\right)$, $n \in \mathbb{N}_+$, with $T_{\theta}^0(x)=x$
and $\vartheta_1=\vartheta_1(x)=\left\lfloor \frac{1}{x \theta} \right\rfloor$ are called the \textit{digits} of $x$ with respect to the $\theta-$expansion in (\ref{3.10}), and we have that the finite truncation of (\ref{3.10}),
$p_{\theta,n}/q_{\theta,n}:=[\vartheta_1, \vartheta_2, \ldots, \vartheta_n]_{\theta}$, tends to $x$ as $n \rightarrow \infty$.

If $\theta^2 = 1/{s}$,  $s \in \mathbb{N}_+$, the digits $\vartheta_n$'s take values greater or equal to $s$, and the transformation $T_{\theta}$ is ergodic with respect to an absolutely continuous invariant probability measure
\begin{equation} \label{3.11}
G_{\theta} (A) :=
\frac{1}{\log \left(1+\theta^{2}\right)} \int_{A} \frac{\theta \mathrm{d}x}{1 + x \theta}, \quad A \in {\mathcal{B}}_{[0, 1]}.
\end{equation}

An $n$-block $(\vartheta_1, \vartheta_2, \ldots, \vartheta_n)$ is said to be \textit{admissible} for the expansion in (\ref{3.10}) if there exists
$x \in [0, \theta)$ such that $\vartheta_i(x)=\vartheta_i$ for all $1 \leq i \leq n$.
If $(\vartheta_1, \vartheta_2, \ldots, \vartheta_n)$ is an admissible sequence, we call the set
\begin{equation} \label{3.013}
I_{\theta} (\vartheta_1, \vartheta_2, \ldots, \vartheta_n) = \{x \in [0, \theta]:  \vartheta_1(x) = \vartheta_1, \vartheta_2(x) = \vartheta_2, \ldots, \vartheta_n(x) = \vartheta_n \},
\end{equation}
\textit{the $n$-th order cylinder}.

Define $(u_{\theta,i})_{i \geq s}$ by
\begin{equation}\label{3.014}
 u_{\theta,i}: [0, \theta] \rightarrow [0, \theta]; \quad   u_{\theta,i}(x) := \frac{1}{i \theta+x}.
\end{equation}
For each $i \geq s$, $u_{\theta,i}=\left( \left.T_{\theta}\right|_{I_{\theta}(a_1)} \right)^{-1}$.
Let
\begin{equation}\label{3.015}
 u_{\theta,\vartheta_1\vartheta_2 \ldots \vartheta_n}(t) := \left(u_{\theta,\vartheta_1} \circ u_{\theta,\vartheta_2} \circ \ldots \circ u_{\theta,\vartheta_n}\right)(t) =
 \frac{1}{\displaystyle \vartheta_1 \theta +\frac{1}{\displaystyle \vartheta_2 \theta + \ddots +\frac{1}{\vartheta_n \theta+t}}}.
\end{equation}
We observe that $u_{\theta,\vartheta_1\vartheta_2 \ldots \vartheta_n}=\left( \left.T^n_{\theta}\right|_{I_{\theta}(\vartheta_1, \vartheta_2, \ldots, \vartheta_n)} \right)^{-1}$.
Therefore
\begin{equation} \label{3.016}
I_{\theta} (\vartheta_1, \vartheta_2, \ldots, \vartheta_n) = \{ u_{\theta,\vartheta_1 \ldots \vartheta_n}(t) : t \in [0, \theta) \}
\end{equation}
which is an interval with the endpoints
$\frac{p_{\theta,n}}{q_{\theta,n}}$ and $\frac{p_{\theta,n}+ \theta p_{\theta,n-1}}{q_{\theta,n}+ \theta q_{\theta,n-1}}$.
Such intervals form a partition of $[0, \theta]$.

We now check R\'enyi condition (\ref{1.05}).
We will use directly, without mentioning them here, some properties proved in \cite{SebeLascu-2014}.
Thus, we have
\begin{equation}\label{3.017}
\frac{\left| u'_{\theta,\vartheta_1 \ldots \vartheta_n}(t) \right|}{\left| u'_{\theta,\vartheta_1 \ldots \vartheta_n}(r) \right|} =
\left(\frac{q_{\theta,n}+r q_{\theta,n-1}}{ q_{\theta,n}+t q_{\theta,n-1} }\right)^2 \leq
\left(\frac{q_{\theta,n}+\theta q_{\theta,n-1}}{ q_{\theta,n}}\right)^2 \leq \left( 1+\theta^2 \right)^2.
\end{equation}

We compute the entropy $h(T_{\theta})$ by the Rohlin's formula (\ref{1.1}):
\begin{eqnarray} \label{3.12}
h(T_{\theta}) &=& \int^{\theta}_{0} \log\left|T'_{\theta} (x)\right| \mathrm{d} G_{\theta}(x)
              = \int^{\theta}_{0} \frac{- \log x^2}{\log(1+\theta^2)}\frac{\theta \mathrm{d}x}{1+\theta x} \nonumber \\
              &=&  \frac{-2 \theta}{\log(1 + \theta^2)} \int_{0}^{\theta} \frac{\log x}{1 + \theta x}\mathrm{d}x.
\end{eqnarray}
As examples, we have the values of $h(T_{\theta})$ for different values of $s=1/\theta^2$ and the graph of $h(T_{\theta})$ in Table~\ref{tab1}.

\subsection{$N-$continued fraction expansions} \label{sec3.3}
Fix an integer $N \geq 1$.
The measure-theoretical dynamical system $([0, 1],{\mathcal B}_{[0, 1]}, T_N, G_N)$ is defined as follows:
\begin{equation} \label{3.13}
T_N:[0, 1] \rightarrow [0, 1]; \quad
T_{N}(x):=
\left\{
\begin{array}{ll}
{\displaystyle \frac{N}{x}- \left\lfloor\frac{N}{x}\right\rfloor},&
{ \mbox{if } x \in (0, 1],}\\
0,& \mbox{if }x=0
\end{array}
\right.
\end{equation}
and
\begin{equation} \label{3.14}
G_N (A) :=
\frac{1}{\log \frac{N+1}{N}} \int_{A} \frac{\mathrm{d}x}{x+N},
\quad A \in {\mathcal{B}}_{[0, 1]}.
\end{equation}
The probability measure $G_N$ is $T_N$-invariant, and the dynamical system \linebreak
$([0, 1],{\mathcal B}_{[0, 1]}, T_N, G_N)$ is ergodic.

For any $0<x<1$ put $\varepsilon_1(x) = \lfloor N/x \rfloor$ and $\varepsilon_n(x)=\varepsilon_1\left(T_N^{n-1}(x)\right)$, $n \in \mathbb{N}_+$, with $T_N^{0}(x)=x$.
Then every irrational $0<x<1$ can be written in the form
\begin{equation} \label{3.15}
x = \displaystyle \frac{N}{\varepsilon_1+\displaystyle \frac{N}{\varepsilon_2+\displaystyle \frac{N}{\varepsilon_3+ \ddots}}} =: [\varepsilon_1, \varepsilon_2, \varepsilon_3, \ldots]_N,
\end{equation}
where $\varepsilon_n$'s are non-negative integers, $\varepsilon_n \geq N$, $n \in \mathbb{N}_+$.
We call (\ref{3.15}) the \textit{$N-$continued fraction expansion} of $x$ and
$p_{N,n}(x) / q_{N,n}(x) := [\varepsilon_1, \varepsilon_2, \ldots, \varepsilon_n]_N$ the \textit{$n$-th order convergent} of $x \in [0, 1]$.
Then $p_{N,n}(x) / q_{N,n}(x) \rightarrow x$, $n \rightarrow \infty$.

An $n$-block $(\varepsilon_1, \varepsilon_2, \ldots, \varepsilon_n)$ is said to be \textit{admissible} for the expansion in (\ref{3.15}) if there exists $x \in [0, 1)$ such that $\varepsilon_i(x)=\varepsilon_i$ for all $1 \leq i \leq n$.
If $(\varepsilon_1, \varepsilon_2, \ldots, \varepsilon_n)$ is an admissible sequence, we call the set
\begin{equation} \label{3.018}
I_N (\varepsilon_1, \varepsilon_2, \ldots, \varepsilon_n) = \{x \in [0, 1]:  \varepsilon_1(x) = \varepsilon_1, \varepsilon_2(x) = \varepsilon_2, \ldots, \varepsilon_n(x) = \varepsilon_n \},
\end{equation}
\textit{the $n$-th order cylinder}.

Define $(u_{N,i})_{i \geq N}$ by
\begin{equation}\label{3.019}
 u_{N,i}: [0, 1] \rightarrow [0, 1]; \quad   u_{N,i}(x) := \frac{N}{i + x}.
\end{equation}
For each $i \geq N$, $u_{N,i}=\left( \left.T_{N}\right|_{I_{N}(\varepsilon_1)} \right)^{-1}$.
Let
\begin{equation}\label{3.020}
 u_{N,\varepsilon_1\varepsilon_2 \ldots \varepsilon_n}(t) :=
 \left(u_{N,\varepsilon_1} \circ u_{N,\varepsilon_2} \circ \ldots \circ u_{N,\varepsilon_n}\right)(t) =
 \frac{N}{\displaystyle \varepsilon_1 +\frac{N}{\displaystyle \varepsilon_2 + \ddots +\frac{N}{\varepsilon_n +t}}}.
\end{equation}
We observe that $u_{N,\varepsilon_1\varepsilon_2 \ldots \varepsilon_n}=\left( \left.T^n_{N}\right|_{I_{N}(\varepsilon_1, \varepsilon_2, \ldots, \varepsilon_n)} \right)^{-1}$.
Therefore
\begin{equation} \label{3.021}
I_{N} (\varepsilon_1, \varepsilon_2, \ldots, \varepsilon_n) = \{ u_{N,\varepsilon_1 \ldots \varepsilon_n}(t) : t \in [0, 1] \}
\end{equation}
which is an interval with the endpoints
$\frac{p_{N,n}}{q_{N,n}}$ and $\frac{p_{N,n} +p_{N,n-1} }{q_{N,n} +q_{N,n-1} }$.
Such intervals form a partition of $[0, 1]$.

We now check R\'enyi condition (\ref{1.05}).
We will use directly, without mentioning them here, some properties proved in \cite{SL-2020-4}.
Thus, we have
\begin{equation}\label{3.022}
\frac{\left| u'_{N,\varepsilon_1 \ldots \varepsilon_n}(t) \right|}{\left| u'_{N,\varepsilon_1 \ldots \varepsilon_n}(r) \right|} =
\left(\frac{q_{N,n}+r q_{N,n-1}}{ q_{N,n}+t q_{N,n-1} }\right)^2 \leq
\left(\frac{q_{N,n} + q_{N,n-1}}{ q_{N,n}}\right)^2 \leq \left( \frac{N+1}{N} \right)^2.
\end{equation}
Using Rohlin's entropy formula (\ref{1.1}), we have:
\begin{eqnarray} \label{3.16}
h(T_{N}) &=& \int^{1}_{0} \log\left|T'_{N} (x)\right| \mathrm{d} G_{N}(x) = \frac{1}{\log \frac{N+1}{N}} \int^{1}_{0} \frac{\log\left(\frac{N}{x^2} \right)}{x+N}\mathrm{d}x \nonumber \\
         &=& \frac{\frac{\pi^2}{3}+2\mathrm{Li}_2(N+1)+\log(N+1)\log N}{{\log \frac{N+1}{N}}},
\end{eqnarray}
where $\mathrm{Li}_2$ denotes the \textit{dilogarithm function}, defined by
\begin{equation}\label{3.17}
  \mathrm{Li}_2(x) = \int_{0}^{x} \frac{\ln t }{1-t} \mathrm{d}t \quad \mbox{ or } \quad \mathrm{Li}_2(x) = \sum_{k=1}^{\infty} \frac{x^k}{k^2}.
\end{equation}
As examples, we have the values of $h(T_{N})$ for different values of $N$ and the graph of $h(T_{N})$ in Table~\ref{tab1}.

\subsection{R\'enyi-type continued fraction expansions} \label{sec3.4}
Fix an integer $N \geq 2$. Let the \textit{R\'enyi-type continued fraction transformation} $R_N : [0, 1] \rightarrow [0, 1]$ be given by
\begin{equation} \label{3.19}
R_{N}(x):=
\left\{
\begin{array}{lll}
{\displaystyle \frac{N}{1-x}- \left\lfloor\frac{N}{1-x}\right\rfloor}, & { \mbox{if } x \in [0, 1),} \\
 0, & \mbox{if } x=1.
\end{array}
\right.
\end{equation}
For any irrational $x \in [0, 1]$, $R_N$ generates a new continued fraction expansion of $x$ of the form
\begin{equation} \label{3.20}
x = 1 - \displaystyle \frac{N}{1+r_1 - \displaystyle \frac{N}{1+r_2 - \displaystyle \frac{N}{1+r_3 - \ddots}}} =:[r_1, r_2, r_3, \ldots]_R.
\end{equation}
Here, $r_n$'s are non-negative integers greater than or equal to $N$ defined by
\begin{equation} \label{3.21}
r_1:=r_1(x) = \left\lfloor \frac{N}{1-x} \right\rfloor, x \neq 1; \quad r_1(1)=\infty
\end{equation}
and
\begin{equation}
r_n := r_n(x) = r_1\left( R^{n-1}_N (x) \right), \quad n \geq 2, \label{3.22}
\end{equation}
with $R_{N}^0 (x) = x$.
The sequence of rationals $\left\{{p_{R,n}}/{q_{R,n}}\right\}:=[r_1, r_2, \ldots, r_n]_R$, ${n \in \mathbb{N_+}}$ are the convergents to $x$ in $[0, 1]$.

The dynamical system $\left([0,1], {\mathcal B}_{[0,1]}, R_N, \rho_N \right)$ is measure preserving and ergodic, where the probability measure $\rho_N$ is defined by
\begin{equation} \label{3.23}
\rho_N (A) :=
\frac{1}{\log \left(\frac{N}{N-1}\right)} \int_{A} \frac{\mathrm{d}x}{x+N-1}, \quad A \in {\mathcal{B}}_{[0,1]}.
\end{equation}
An $n$-block $(r_1, r_2, \ldots, r_n)$ is said to be \textit{admissible} for the expansion in (\ref{3.20}) if there exists $x \in [0, 1)$ such that $r_i(x)=r_i$ for all $1 \leq i \leq n$.
If $(r_1, r_2, \ldots, r_n)$ is an admissible sequence, we call the set
\begin{equation} \label{3.026}
I_R (r_1, r_2, \ldots, r_n) = \{x \in [0, 1]:  r_1(x) = r_1, r_2(x) = r_2, \ldots, r_n(x) = r_n \},
\end{equation}
\textit{the $n$-th order cylinder}.

Define $(u_{R,i})_{i \geq N}$ by
\begin{equation}\label{3.038}
 u_{R,i}: [0, 1] \rightarrow [0, 1]; \quad   u_{R,i}(x) := 1 - \frac{N}{i + x}.
\end{equation}
For each $i \geq N$, $u_{R,i}=\left( \left.R_{N}\right|_{I_{R}(r_1)} \right)^{-1}$.
Let
\begin{equation}\label{3.039}
 u_{R,r_1r_2 \ldots r_n}(t) :=
 \left(u_{R,r_1} \circ u_{R,r_2} \circ \ldots \circ u_{R,r_n}\right)(t) =
 1 - \frac{N}{\displaystyle 1+r_1 -\frac{N}{\displaystyle 1+ r_2 - \ddots -\frac{N}{r_n +t}}}.
\end{equation}
We observe that $u_{R,r_1r_2 \ldots r_n}=\left( \left.R^n_{N}\right|_{I_{R}(r_1, r_2, \ldots, r_n)} \right)^{-1}$.
Therefore
\begin{equation} \label{3.040}
I_{R} (r_1, r_2, \ldots, r_n) = \{ u_{R,r_1 \ldots r_n}(t) : t \in [0, 1] \}
\end{equation}
which is the interval $\left[ \frac{p_{R,n}-p_{R,n-1}}{q_{R,n}-q_{R,n-1}}, \frac{p_{R,n}}{q_{R,n}} \right)$.
Such intervals form a partition of $[0, 1]$.

Before applying Rohlin's formula, we must check R\'enyi condition (\ref{1.05}).
We will use directly, without mentioning them here, some properties proved in \cite{LS-2020-1}.
Thus, we have
\begin{equation}\label{3.041}
\frac{\left| u'_{R,r_1 \ldots r_n}(t) \right|}{\left| u'_{R,r_1 \ldots r_n}(r) \right|} =
\left(\frac{q_{R,n}+(r-1) q_{R,n-1}}{ q_{R,n}+(t-1) q_{R,n-1} }\right)^2 \leq
\left(\frac{q_{R,n}}{ q_{R,n} - q_{R,n-1}}\right)^2 \leq \left(\frac{N}{N-1}\right)^2.
\end{equation}
The entropy $h(R_N)$ is given by
\begin{eqnarray} \label{3.24}
h(R_{N}) &=& \int^{1}_{0} \log\left|R'_{N} (x)\right| \mathrm{d} \rho_{N}(x) = \frac{1}{\log \frac{N}{N-1}} \int^{1}_{0} \frac{\log\frac{N}{(1-x)^2}}{x+N-1}\mathrm{d}x \nonumber \\
         &=& \log N + \frac{2 \mathrm{Li}_2\left(\frac{1}{N}\right)}{\log \frac{N}{N-1}},
\end{eqnarray}
where $\mathrm{Li}_2$ is as in (\ref{3.17}).
As examples, we have the values of $h(R_{N})$ for different values of ${N}$ and the graph of $h(R_{N})$ in Table~\ref{tab1}.

\begin{table}[H]
\caption{Entropies $h(T_{\ell})$, $h(T_{\theta})$, $h(T_{N})$ and $h(R_{N})$ for different values of the parameters involved} \label{tab1}
\begin{tabular}{cc|cc|cc|cc}
  \hline
  $\quad \ell $ & $ \quad h(T_{\ell})$ & $\quad s $ & $ \quad h(T_{\theta})$ & $\quad N $ & $ \quad h(T_{N})$ & $\quad N $ & $ \quad h(R_{N})$ \\
  \hline
  $\ell=2$   & $1.62258 $   & $s=1$   & $2.37314$  & $N=1$   & $2.37314$  & $N=2$   & $2.37314$  \\
  $\ell=3$   & $1.26775 $   & $s=3$   & $3.24705$  & $N=3$   & $3.24705$  & $N=3$   & $2.905$    \\
  $\ell=5$   & $0.996315 $  & $s=5$   & $3.70244$  & $N=5$   & $3.70244$  & $N=5$   & $3.50063$  \\
  $\ell=10$  & $0.765943 $  & $s=10$  & $4.35074$  & $N=10$  & $4.35074$  & $N=10$  & $4.25052$  \\
  $\ell=50$  & $0.476521 $  & $s=50$  & $5.92195$  & $N=50$  & $5.92195$  & $N=50$  & $5.90194$  \\
  $\ell=100$ & $0.406218 $  & $s=100$ & $6.61015$  & $N=100$ & $6.61015$  & $N=100$ & $6.60015$  \\
  $\ell=200$ & $0.350849 $  & $s=1000$ & $8.90825$ & $N=1000$ & $8.90825$ & $N=1000$ & $8.90726$ \\
  \hline
\end{tabular}
\end{table}


\section{Comparing the efficiency of some expansions} \label{sec4}
In this section we apply the Extended Lochs' theorem presented in Section \ref{sec2.4}
and compare two by two the expansions presented in the previous section.
First of all, we observe that for various values of the parameters involved, the entropies $h(T_{\theta})$ and $h(T_{N})$ are equal.
Since entropy is an isomorphism invariant, we conjecture the following result.

\begin{conjecture}
For an irrational $\theta \in (0, 1)$ and a non-negative integer $N \geq 2$ with $1/{\theta^2}=N$, the transformations $T_{\theta}$ in (\ref{3.9}) and $T_{N}$ in (\ref{3.13}) are isomorphic.
\end{conjecture}

For this reason, we make only the following pairs: $N-$continued fractions and Chan's continued fractions, $N-$continued fractions and R\'enyi-type continued fractions, R\'enyi-type continued fractions and Chan's continued fractions.

We observe that the transformations $T_N$, $T_{\ell}$ and $R_N$ satisfy the two conditions from Extended Lochs' theorem (see Section \ref{sec2.4}).

\subsection{$N-$continued fractions and Chan's continued fractions} \label{sec4.1}
Let $I^x_N (\varepsilon_1, \varepsilon_2, \ldots, \varepsilon_n )$ denotes the $n$-th order cylinder of the $N-$continued fraction that contains $x$, and $I^x_{\ell} (a_1, a_2, \ldots, a_m )$ denotes the $m$-th order cylinder of the Chan's continued fraction that contains $x$.
Then
\begin{equation}\label{4.1}
  m_{N{\ell}} (n,x) := sup\left\{ m:  I^x_N (\varepsilon_1, \varepsilon_2, \ldots, \varepsilon_n ) \subset I^x_{\ell} (a_1, a_2, \ldots, a_m ) \right\}
\end{equation}
represents the number of digits in the Chan's expansion of $x$ in (\ref{3.1}) that can be determined from knowing the first $n$ digits in the $N-$continued fraction in (\ref{3.15}).
Therefore, applying (\ref{2.15}) we have
\begin{equation} \label{4.2}
\lim_{n \rightarrow \infty} \frac{m_{N{\ell}}(n, x)}{n} = \frac{h(T_N)}{h(T_{\ell})}
\end{equation}
where $h(T_N)$ and $h(T_{\ell})$ are as in (\ref{3.16}) and (\ref{3.8}), respectively.
Given the values in Table~\ref{tab1}, we observe that $N$-continued fraction expansion is more effective
than Chan's continued fraction expansion regardless of the values taken by the parameters $N$ and $\ell$, respectively.

As examples, we have,
\begin{equation} \label{4.03}
\lim_{n \rightarrow \infty} \frac{m_{1{2}}(n, x)}{n} = 1.462571953\ldots \quad \mbox{or} \quad
\lim_{n \rightarrow \infty} \frac{m_{3{2}}(n, x)}{n} = 2.001164812\ldots
\end{equation}
So roughly, if we approximate a number from the unit interval by keeping the first $1000$ digits of the $N$-continued fraction expansion, in order to retain the same degree of accuracy we need to keep about $1462$ digits in the Chan's continued fraction expansion.

\subsection{$N-$continued fractions and R\'enyi-type continued fractions} \label{sec4.2}
Let $I^x_N (\varepsilon_1, \varepsilon_2, \ldots, \varepsilon_n )$ denotes the $n$-th order cylinder of the $N-$continued fraction that contains $x$, and $I^x_{R} (r_1, r_2, \ldots, r_m )$ denotes the $m$-th order cylinder of the R\'enyi-type continued fraction that contains $x$.
Then
\begin{equation}\label{4.3}
  m_{NR} (n,x) := sup\left\{ m:  I^x_N (\varepsilon_1, \varepsilon_2, \ldots, \varepsilon_n ) \subset I^x_{R} (r_1, r_2, \ldots, r_m ) \right\}
\end{equation}
represents the number of digits in the R\'enyi-type continued fraction of $x$ in (\ref{3.20}) that can be determined from knowing the first $n$ digits in the $N-$continued fraction in (\ref{3.15}).
Therefore, applying (\ref{2.15}) we have
\begin{equation} \label{4.4}
\lim_{n \rightarrow \infty} \frac{m_{NR}(n, x)}{n} = \frac{h(T_N)}{h(R_N)}
\end{equation}
where $h(T_N)$ and $h(R_N)$ are as in (\ref{3.16}) and (\ref{3.24}), respectively.
Given the values in Table~\ref{tab1}, we observe that $N-$continued fraction expansion is more effective
than R\'enyi-type continued fraction expansion regardless of the values taken by the parameter $N$.

We notice that $h(T_1) = h(R_2)$.
We also have
\begin{equation} \label{4.04}
\lim_{n \rightarrow \infty} \frac{m_{3{3}}(n, x)}{n} = 1.117745267\ldots \quad \mbox{or} \quad
\lim_{n \rightarrow \infty} \frac{m_{5{5}}(n, x)}{n} = 1.057649623\ldots
\end{equation}
As $N$ grows, the entropies are very close, which means that the efficiency of the two continued fraction expansions are about the same.

\subsection{R\'enyi-type continued fractions and Chan's continued fractions} \label{sec4.3}
Let $I^x_{R} (r_1, r_2, \ldots, r_n )$ denotes the $n$-th order cylinder of the R\'enyi-type continued fraction that contains $x$, and
$I^x_{\ell} (a_1, a_2, \ldots, a_m )$ denotes the $m$-th order cylinder of the Chan's continued fraction that contains $x$.
Then
\begin{equation}\label{4.5}
  m_{R{\ell}} (n,x) := sup\left\{ m:  I^x_{R} (r_1, r_2, \ldots, r_n ) \subset I^x_{\ell} (a_1, a_2, \ldots, a_m ) \right\}
\end{equation}
represents the number of digits in the Chan's continued fraction of $x$ in (\ref{3.1}) that can be determined from knowing the first $n$ digits in the R\'enyi-type continued fraction in (\ref{3.20}).
Therefore, applying (\ref{2.15}) we have
\begin{equation} \label{4.6}
\lim_{n \rightarrow \infty} \frac{m_{R{\ell}}(n, x)}{n} = \frac{h(R_N)}{h(T_{\ell})}
\end{equation}
where $h(R_N)$ and $h(T_{\ell})$ are as in (\ref{3.24}) and (\ref{3.8}), respectively.
Given the values in Table~\ref{tab1}, we observe that R\'enyi-type continued fraction expansion is more effective
than Chan's continued fraction expansion regardless of the values taken by the parameters $N$ and $\ell$, respectively.

As examples, we have,
\begin{equation} \label{4.06}
\lim_{n \rightarrow \infty} \frac{m_{2{2}}(n, x)}{n} = 1.462571953\ldots \quad \mbox{or} \quad
\lim_{n \rightarrow \infty} \frac{m_{2{3}}(n, x)}{n} = 1.871930586\ldots
\end{equation}

\section{Final remarks} \label{sec5}
To conclude, $N$-continued fractions are more efficient than R\'enyi-type continued fractions at representing a number in the unit interval.
Since the entropies $h(T_N) \geq 2.37314$ for $N \geq 1$, $h(R_N) \geq 2.37314$ for $N \geq 2$ and $h(T_G) = \pi^2 /(6 \log2)= 2.37314$,
it follows that $N$-continued fractions and R\'enyi-type continued fractions are more efficient than regular continued fractions (RCFs).
Since the entropy $h(T_{\ell}) \leq 1.62258$ for $\ell \geq 2$, it follows that RCFs are more efficient than Chan's continued fractions.
Thus, $N$-continued fractions are the most efficient at representing a number in the unit interval, with a very close efficiency being R\'enyi-type continued fractions.

Our paper is a systematic presentation of continued fraction expansions that have been investigated by us during the last ten years.
Having always in view the classical RCF expansion, we payed attention to some interval maps that generate expansions and that admit an invariant density with suitable ergodic properties. However these ergodic properties, recently studied, are not enough to yield rates of convergence for mixing properties. For this a Gauss-Kuzmin-type theorem is needed. There are still many open questions closely related to this problem.
On the other hand, to extend our references list \cite{Lascu-2013,Sebe-2010,SebeLascu-2014,Sebe-2017,Lascu-2016,SL-2020-4,LS-2020-1,LS-2020-2,SL-2020-3}, we considered the opportunity of starting new investigations such that the efficiency of pairs of maps for which the generating partition has finite entropy. In these circumstances, we reviewed sufficient conditions on the latter to belong of a class of number-theoretic fibred maps for which the generating partition has finite entropy.


\end{document}